\documentclass[12pt]{amsart}

\usepackage[centertags]{amsmath}
\usepackage{amsfonts,amssymb,amsthm}
\usepackage{graphicx}

\theoremstyle{plain}
\newtheorem{thm}{Theorem}[section]
\newtheorem{lem}[thm]{Lemma}
\newtheorem{prop}[thm]{Proposition}
\newtheorem{cor}[thm]{Corollary}

\theoremstyle{definition}
\newtheorem{defn}[thm]{Definition}
\newtheorem{rem}[thm]{Remark}

\begin{document}
\title{\large Upper Bounds on the Number of Vertices of Weight $\leq k$ in Particular Arrangements of Pseudocircles}
\author{\bf Ronald ORTNER}
\email{rortner@unileoben.ac.at}
\address{\parbox{1.4\linewidth}{University of Leoben\\
Franz-Joseph-Strasse 18, A-8700 Leoben, Austria}}
\date{June 13, 2008}

\begin{abstract}
In arrangements of pseudocircles (Jordan curves) the weight of a vertex (intersection point) is the number
of pseudocircles that contain the vertex in its interior.
We give improved upper bounds on the number of vertices of weight $\leq k$ in certain arrangements of pseudocircles
in the plane. In particular, forbidding certain subarrangements we improve the known bound of $6n-12$ (cf.\ \cite{kede})
for vertices of weight 0 in arrangements of $n$ pseudocircles to $4n-6$. In \textit{complete} arrangements (i.e.\ 
arrangements with each two pseudocircles intersecting) we identify two subarrangements of three and four pseudocircles, respectively, whose absence gives improved bounds for vertices of weight 0 and more generally for vertices
of weight $\leq k$.
\end{abstract} 
\maketitle

\markboth{\sf R. Ortner}{\sf Bounds on the Number of Vertices of Weight $\leq k$ in Arrangements of Pseudocircles}
\baselineskip 15pt

\section{Introduction}
A \emph{pseudocircle} is a simple closed (Jordan) curve in the plane.
An \emph{arrangement of pseudocircles} is a finite set $\Gamma=\{\gamma_1, \ldots, \gamma_n\}$ of simple closed 
curves in the plane such that
\begin{itemize}
  \item[$\rhd$] no three curves meet each other at the same point,
\item[$\rhd$] each two curves $\gamma_i, \gamma_j$ have at most two points in common, and
\item[$\rhd$]  these \textit{intersection points} in $\gamma_i\cap\gamma_j$ are always 
points where $\gamma_i, \gamma_j$ cross each other.
\end{itemize}  
An arrangement is \emph{complete} if each two pseudocircles intersect. 

\smallskip

Any arrangement can be interpreted as a planar embedding of a graph
whose vertices are the intersection points between the pseudocircles and whose
edges are the curves between these intersections. In the following we will often
refer to this graph when talking about \textit{vertices}, \textit{edges}, and 
\textit{faces} of the arrangement.

\begin{defn}
Let $\Gamma=\{\gamma_1, \ldots, \gamma_n\}$ be an arrangement of pseudocircles. 
The \emph{weight of a vertex} $V$ is the number of pseudocircles
$\gamma_i$ such that $V$ is contained in ${\rm int}(\gamma_i)$, the interior  of
$\gamma_i$. Weights of edges and faces are defined accordingly.

We will consider the  number $v_k=v_k(\Gamma)$ of vertices of given weight $k$,
the number $v_{\leq k}=v_{\leq k}(\Gamma)$ of vertices of weight $\leq k$, and
the number $v_{\geq k}=v_{\geq k}(\Gamma)$ of vertices of weight $\geq k$.
Further, $f_k=f_k(\Gamma)$ denotes the number of faces of weight $k$. 
\end{defn}

Concerning the characterization of the \textit{weight vectors} $(v_0,v_1,\ldots,$ $v_{n-2})$ 
of arrangements of pseudocircles little is known.
So far, sharp upper bounds on the $v_k$ exist only for $k=0$.

\begin{thm}[Kedem et al. \cite{kede}]\label{thm:0}
For all arrangements $\Gamma$ with $n:=|\Gamma| \geq 3$:
\[ v_0 \;\leq\; 6n-12.\]
Moreover, for each $n \geq 3$ there is an arrangement of $n$ (proper) circles in the plane such that $v_0 = 6n-12.$
\end{thm}

Theorem \ref{thm:0} can be used to obtain general upper bounds on $v_{\leq k}$ by
some clever probabilistic methods.

\begin{thm}[Sharir \cite{shar}]\label{thm:<k}
For all arrangements of $n$ pseudocircles and all $k > 0$:
\[ v_{\leq k}\; \leq\; 26kn.\]
\end{thm}

On the other hand, J.\ Linhart and Y.\ Yang established the following sharp upper bound on $v_{\geq k}$.

\begin{thm}[Linhart, Yang \cite{linh}]\label{thm:>k}
For all arrangements of $n \geq 2$ pseudocircles and all $k$ with $0 \leq k \leq n-2$:
\[ v_{\geq k}\; \leq\; (n+k)(n-k-1).\]
\end{thm}

Beyond Theorems \ref{thm:0} -- \ref{thm:>k} there are hardly any results for arrangements of pseudocircles in general.
 In this paper we are going to improve the upper bounds of Theorems \ref{thm:0} and \ref{thm:<k} for some particular classes of 
arrangements.

\smallskip

\section{Preparations and Some Minor Results}
In this section we show some minor results based on Theorems \ref{thm:>k} and \ref{thm:0}. 
\subsection{Improved Bounds by Theorem \ref{thm:>k}}
We start with a bound due to J.\ Linhart and based upon a result of Y.\ Yang which holds, 
if there is a face of large weight in the arrangement.

\begin{prop}[Yang \cite{yang}]\label{prop:inout}
Let $\Gamma$ be an arrangement of $n$ pseudocircles in the plane with weight vector $(v_0,v_1, \ldots, v_{n-2})$ 
and $f_{n}>0$. Then there is an arrangement of pseudocircles $\Gamma'$ with weight vector
$(v'_0,v'_1, \ldots, v'_{n-2})=(v_{n-2},v_{n-1}, \ldots, v_0)$.
\end{prop}

\begin{proof}
Choose an arbitrary point $P$ contained in the face of weight $n$. We place a sphere on the plane touching it in $P$. By stereographic projection from the antipodal point of $P$, $P^*$, we obtain an arrangement $\Gamma_{S^2}$ of pseudocircles on $S^2$. Obviously, when projecting an arrangement from the plane to the sphere, we want the interiors of pseudocircles in the plane to be projected to the interiors of pseudocircles on the sphere. Thus, the arrangement on the sphere has the same weight vector as $\Gamma$. Now consider the arrangement $\Gamma'_{S^2}$ arising when we swap interior and exterior of each pseudocircle. Then an arbitrary point $Q$ on $S^2$ is contained in the interior of $\gamma$ in $\Gamma'_{S^2}$ if and only if $Q$ is outside $\gamma$ in $\Gamma_{S^2}$. Therefore a vertex of weight $k$ in $\Gamma'_{S^2}$ has weight $n-k-2$ in $\Gamma_{S^2}$ and vice versa. In particular, $P$ is not contained in the interior of any pseudocircle in $\Gamma'_{S^2}$. Thus, we can project the whole arrangement from $P$ to the plane tangent to $P^*$ yielding an arrangement $\Gamma'$ with the claimed property.
\end{proof}

J.\ Linhart pointed out to me that Proposition \ref{prop:inout} together with Theorem \ref{thm:>k} yields the following improvement of the upper bound on $v_{\leq k}$ for arrangements with $f_n>0$.

\begin{thm}[Linhart \cite{linh2}]\label{thm:fn}
For all arrangements $\Gamma$ of $n$ pseudocircles with $f_n >0$:
\[ v_{\leq k} \; \leq \; 2(k+1)n - (k+1)(k+2). \]
\end{thm}

\begin{proof}
Let $\Gamma$ be an arrangement with weight vector $(v_0,v_1, \ldots, v_{n-2})$ and $f_n>0$. 
Then by Proposition \ref{prop:inout}, there exists an arrangement $\Gamma'$ with $v'_k=v_{n-k-2}$ vertices of 
weight $k$ for $0 \leq k \leq n-2$. Therefore,
\[ v_{\leq k} \; = \; \sum_{j=0}^k v_j \; = \; \sum_{j=0}^k v'_{n-j-2} \; =  \sum_{j=n-k-2}^{n-2} v'_j \; = \; v'_{\geq n-k-2}. \]
Applying Theorem \ref{thm:>k} yields
\begin{eqnarray*}
	 v_{\leq k} &=& v'_{\geq n-k-2} \; \leq \; (n+n-k-2)(n-(n-k-2)-1)\; =  \\
	 &=& (2n-k-2)(k+1) \; = \; 2(k+1)n - (k+1)(k+2). \hspace{1cm}\qedhere
\end{eqnarray*}
\end{proof}

\subsection{Improved Bounds by Theorem \ref{thm:0}}
\quad On the other hand, bounds on $v_0$ can also be improved if there is a face of weight 0
with many pseudocircles participating in its boundary.

\begin{prop}\label{prop:many}
Let $\Gamma$ be an arrangement of $n$ pseudocircles with a face $F$ of weight 0 
such that for each $\gamma\in\Gamma$ there is an edge of $\gamma$ on $\partial F$, the boundary of $F$. 
Then
\[ v_0 \; \leq \; 4n-6. \]
\end{prop}

\begin{proof}
Let us assume that there exists an arrangement $\Gamma$ as described in the proposition such that $v_0>4n-6$. 
As shown in Figure \ref{fig:construction} we can add a pseudocircle $\gamma'$ to $\Gamma$ 
such that $\gamma'$ cuts each $\gamma\in\Gamma$ on $\partial F$ in two vertices of weight 0.
\begin{figure}[htbp]
	\centering
		\scalebox{0.50}{\includegraphics{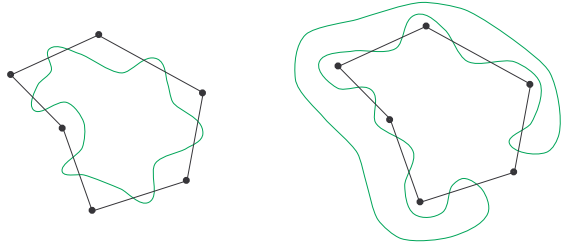}}
	\caption{\label{fig:construction} Adding a pseudocircle $\gamma'$ cutting each $\gamma\in\Gamma$ in two vertices of weight 0.}
\end{figure}
 Note that we may add $\gamma'$ such that it does not contain any vertices of $\Gamma$ in its interior.
 Hence, in the arrangement $\Gamma':=\Gamma\cup\{\gamma'\}$ we have
\[ v_0(\Gamma') \; > \; 4n-6 + 2n\; = \; 6(n+1)-12,\]
which contradicts Theorem \ref{thm:0}.
\end{proof}

\begin{figure}[h]
	\centering
		\scalebox{0.50}{\includegraphics{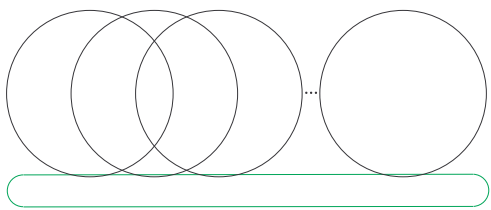}}
	\caption{\label{fig:facemax} Arrangement of $n$ pseudocircles with $v_0=4n-6$.}
\end{figure}

The bound of Proposition \ref{prop:many} is also sharp, see e.g.\ the arrangement shown in Figure \ref{fig:facemax}.

\begin{rem}
The proof method of Proposition \ref{prop:many} can evidently be generalized to give a general bound of
$v_0\leq 6n-2k-6$ for arrangements of $n$ pseudocircles with a face in whose boundary $k$
pseudocircles participate.
\end{rem}

\bigskip

Theorem \ref{thm:0} can also be used to obtain an upper bound on $v_0$ that depends on $f_0$.

\begin{thm}\label{thm:v0f0}
Let $\Gamma$ be an arrangement of $n$ pseudocircles. Then
\[ v_0 \; \leq \; 2n + 2f_0 - 4. \]
\end{thm}

\begin{figure}[b]
\setlength{\unitlength}{1cm}
\begin{picture}(13,2)
  \put(3,1){\circle{2}}
  \put(4,1){\circle{2}}
  \put(5,1){\circle{2}}
  \put(6,1){\circle{2}}
  
  \put(6.75,1){$\ldots$}
  
  \put(8,1){\circle{2}}
  \put(9,1){\circle{2}}
  \put(10,1){\circle{2}}
  \put(11,1){\circle{2}}
\end{picture}
  \caption{\label{fig:f0} Arrangement of $n$ circles with $f_0=1$ and $v_0=2n-2$ (cf.\ Theorem \ref{thm:v0f0}).}
\end{figure}
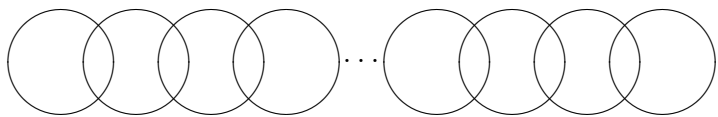

Theorem \ref{thm:v0f0} can be proved without much effort from Theorem \ref{thm:0} with the aid of
the following upper bound on $f_0$ in Proposition \ref{prop:f0}, which is also an easy consequence of Theorem \ref{thm:0}.
As Theorem \ref{thm:v0f0} together with Proposition \ref{prop:f0} is obviously a generalization of
Theorem \ref{thm:0}, this can be considered as self-strengthening of Theorem \ref{thm:0}.

\begin{prop}\label{prop:f0}
Let $\Gamma$ be an arrangement of $n\geq 3$ pseudocircles. Then
\[ f_0 \; \leq \; 2n - 4. \]
\end{prop}
\begin{proof}
First note that the boundary of each face of weight 0 consists of at least three edges (and hence vertices) of weight 0. 
For if there were a face with only two edges belonging to some pseudocircles $\gamma_i$ and 
$\gamma_j$, then $\gamma_i \cap\gamma_j$ would have more than the two allowed intersection points.
On the other hand, each vertex of weight 0 is on the boundary of only a single face of weight 0.
Therefore by Theorem \ref{thm:0},
\[  f_0 \,\leq\, \frac{v_0}{3}\,\leq\, \frac{6n-12}{3}\,=\,2n-4. \qedhere \]
\end{proof}

\begin{proof}[Proof of Theorem \ref{thm:v0f0}]
(\textit{Note:} The induction proof given in prior versions of this paper does not work.)
If all faces of weight 0 are triangles, 
then $v_0=3f_0$, so that by Proposition \ref{prop:f0},
\[ v_0 \;=\; 2 f_0 + f_0 \;\leq\; 2f_0 + 2n-4.  \]
We proceed by induction on the number of faces of weight 0 with $\ell>3$ edges.
Adding a pseudocircle $\gamma$ as described in the proof of Proposition \ref{prop:many}
removes such a face, while adding $\ell-1$ new faces of weight 0 and $2\ell$ new vertices 
of weight 0. Then by induction assumption
\[ v_0 + 2\ell \;\leq\; 2(n+1) + 2(f_0+(\ell-1)) - 4,  \]
whence the theorem follows.
\end{proof}

\bigskip

\section{New Bounds for Complete Arrangements with Forbidden Subarrangements}

\begin{figure*}[htbp]
	\centering
		\scalebox{0.40}{\includegraphics{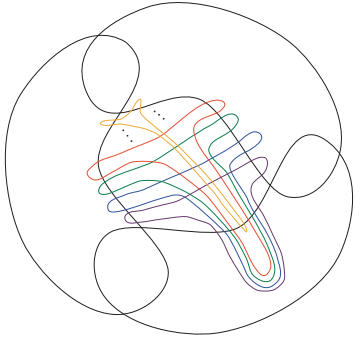}}
	\caption{\label{fig:maxcomplete} A complete arrangement of pseudocircles with $v_0=6n-12$.}
\end{figure*}

In this section we consider bounds on complete arrangements. First, we'd like to remark that the bound of 
Theorem \ref{thm:0} is sharp for complete arrangements, too. 
That is, for each $n\geq 3$ there is a complete arrangement with $v_0 = 6n-12$ 
(see Figure \ref{fig:maxcomplete}). 

Even for arrangements of unit circles no significant improvement is possible. Consider the arrangement 
obtained from a densest circle packing by increasing the radius of the circles by a sufficiently small value,
so that each touching point is transformed into two intersection points. Then each circle intersects six others,
which gives a total of $\approx 6n$ vertices of weight 0.

\subsection{Forbidding $\alpha$-subarrangements}
\begin{figure}[h]
	\centering
		\scalebox{0.70}{\includegraphics{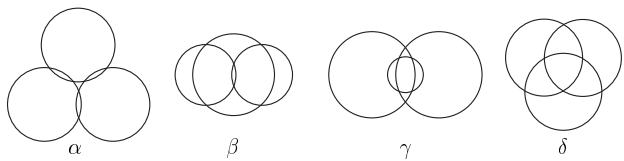}}
	\caption{\label{fig:abcd} Complete arrangements of three pseudocircles in the plane.}\label{fig:3sub}
\end{figure}

Thus, in order to obtain improved bounds on $v_0$, one has to put some additional restrictions
on the arrangement, e.g.\ by forbidding certain subarrangements as we will do it here. 
Evidently, arrangements of three pseudocircles are the smallest subarrangement of interest here.
Figure \ref{fig:abcd} shows the four different  types of subarrangements of three pseudocircles 
one has to take into account. 
Subarrangements of type $\alpha$ play a special role here. Not only are they the only arrangements
of three pseudocircles which meet the bound of Theorem \ref{thm:0}. They are also the only complete 
arrangements of three pseudocircles without any face of weight 3, which is of importance in the light of the 
following Helly type theorem.
The theorem was independently found by E.\ Helly himself and B.\ Ker\'ekj\'art\'o. 
In 1930, Helly published a generalization for dimensions $\geq 2$ (see \cite{hell}).

\begin{thm}[Helly \cite{hell}; Ker\'ekj\'art\'o]\label{thm:helly}
Let $\Gamma=\{\gamma_1,\ldots, \gamma_n \}$ be an
arrangement of pseudocircles such that for all pairwise distinct $\gamma_i,
\gamma_j, \gamma_k:\, {\rm int}(\gamma_i)\cap  {\rm int}(\gamma_j)\cap
 {\rm int}(\gamma_k)\neq \varnothing$. Then
\[ \bigcap_{i=1}^n \, {\rm int}(\gamma_i)\neq \varnothing. \]
\end{thm}

\begin{cor}\label{cor:noalpha}
Let $\Gamma$ be a complete arrangement of $n\geq 2$ pseudocircles that has no subarrangement of type $\alpha$. Then
\[ v_{\leq k} \; \leq \; 2(k+1)n - (k+1)(k+2). \]
\end{cor}

\begin{proof}
Since $\Gamma$ has no $\alpha$-subarrangement, the condition in Theorem \ref{thm:helly} holds, and we may conclude
that there is a face of weight $n$ in $\Gamma$. Applying Theorem \ref{thm:fn} yields the claimed bound.
\end{proof}

\subsection{Forbidding $\alpha^4$-subarrangements}
We have seen that the bounds of Theorems \ref{thm:0} and \ref{thm:<k} can be significantly improved for $\alpha$-free arrangements.
 It is a natural question whether there are alternative bounds for other forbidden subarrangements as well. 
The unique complete arrangement of four pseudocircles that meets the bound of Theorem \ref{thm:0} seems to be
a good candidate. In such an \textit{$\alpha^4$-arrangement} each subarrangement of three pseudocircles
is of type $\alpha$. $\alpha^4$-arrangements prominently appear in the arrangement of Figure \ref{fig:maxcomplete}, 
where the three outer pseudocircles together with any other pseudocircle form an $\alpha^4$-arrangement. 
Indeed, for $\alpha^4$-free arrangements in which there is also no $\beta$-subarrangement 
we can show the following improved upper bound on $v_0$.

\begin{thm}\label{thm:no4alpha}
In complete arrangements of $n\geq 2$ pseudocircles that are $\alpha^4$-free and $\beta$-free, 
\[ v_{0} \; \leq \; 4n-6. \]
\end{thm}

Theorem \ref{thm:no4alpha} follows immediately from the following bound on $f_0$ 
together with Theorem \ref{thm:v0f0}.

\begin{thm}\label{thm:f0c}
In complete arrangements of $n\geq 2$ pseudocircles that are $\alpha^4$-free and $\beta$-free, 
$$f_0\leq n-1.$$
\end{thm}

For the proof of Theorem \ref{thm:f0c} the following lemma is useful. We skip a proof.

\begin{lem}\label{lem:b-free}
 Let $\Gamma$ be a complete, $\beta$-free arrangement. 
Then for each face $F$ of weight 0 in $\Gamma$ there 
is a unique $\alpha$-arrangement $\Gamma_\alpha$ in $\Gamma$ such that
$F$ is the bounded face of weight 0 of $\Gamma_\alpha$.
In particular, each face of weight 0 has only three edges.
\end{lem}

\noindent
\textbf{Proof of Theorem \ref{thm:f0c}.}
We give a proof by induction on $n=|\Gamma|$. The case $n=2$ is trivial, while
for $n=3$ one may consult Figure \ref{fig:3sub}. If $n>3$, choose an
arbitrary pseudocircle $\gamma$ in $\Gamma$. By induction assumption the theorem
holds for $\Gamma':=\Gamma\setminus\{\gamma\}$. We claim that adding $\gamma$
to $\Gamma'$ will increase $f_0$ by at most 1. Indeed, $f_0$ could be increased
by more than 1 only in one of the following two cases: 

First, $\gamma$ may
separate a single face $F$ of weight 0 in $\Gamma'$ into more than two new faces
of weight 0. By Lemma \ref{lem:b-free} such a face $F$ has only three edges
which belong to three pseudocircles that form an $\alpha$-arrangement $\Gamma_\alpha$.
Thus, in order to separate $F$, $\gamma$ has to intersect each pseudocircle of 
$\Gamma_\alpha$ in two vertices of weight 0, so that $\Gamma_\alpha\cup\{\gamma\}$
would be a forbidden $\alpha^4$-arrangement.

On the other hand, there might be two distinct faces $F_1$, $F_2$ in $\Gamma'$,
such that $\gamma$ separates each $F_i$ into two new faces of weight 0.
By Lemma \ref{lem:b-free}, there is a unique $\alpha$-arrangement $\Gamma_\alpha$
enclosing $F_1$, so that $F_2$ will be outside $\Gamma_\alpha$ (i.e. contained in the 
unbounded face of weight 0). Hence, $\gamma$ would have to intersect the bounded
as well as the unbounded face of weight 0 of $\Gamma_\alpha$.
But it is easy to see that this can only happen if $\gamma$ together with two
pseudocircles in $\Gamma_\alpha$ forms a forbidden $\beta$-subarrangement.$\hfill\Box$\\

The bound on $v_0$ of Theorem \ref{thm:no4alpha} is sharp.
Take $n-1$ pseudocircles such that any subarrangement of three pseudocircles
is of type $\delta$. In this arrangement $f_0=1$ and each pseudocircle
has an edge (and hence two vertices) on the single face of weight 0. Adding another pseudocircle
just as inducated in Figure \ref{fig:construction} yields an arrangement with $f_0=n-1$ and 
$v_0=4n-6$.

\smallskip

The improved upper bound on $v_0$ of Theorem \ref{thm:no4alpha}
can be used to improve the upper bound on $v_{\leq k}$ for complete, 
$\alpha^4$-free arrangements.

\begin{thm}\label{thm:k}
For complete, $\alpha^4$-free and $\beta$-free arrangements of $n\geq 2$ pseudocircles and $k>0$,
\[  v_{\leq k} \;\leq\; 18kn.  \]
\end{thm}
\begin{proof}
The proof is basically identical to the proof of Theorem \ref{thm:<k} in \cite{shar},
only with the application of Theorem \ref{thm:0} replaced by an application
of Theorem \ref{thm:no4alpha} and the constants adapted accordingly.
\end{proof}

\begin{rem}
The upper bound of $18kn$ on $v_{\leq k}$ also holds for the class of arrangements mentioned in Proposition
\ref{prop:many}. The proof is again a modification of the original proof in \cite{shar}. Actually, this proof can be adapted
to use any upper bound on $v_0$ for some class of arrangements $\mathcal C$, as long as it is guaranteed that
subarrangements of arrangements in $\mathcal C$ are also in $\mathcal C$. In the proof of Theorem \ref{thm:k}
we use e.g.\ that  subarrangements of $\alpha^4$-free and $\beta$-free arrangements are still $\alpha^4$-free 
and $\beta$-free.
\end{rem}

\section{Conclusion}
\begin{figure}[h]
	\centering
		\scalebox{0.2}{\includegraphics{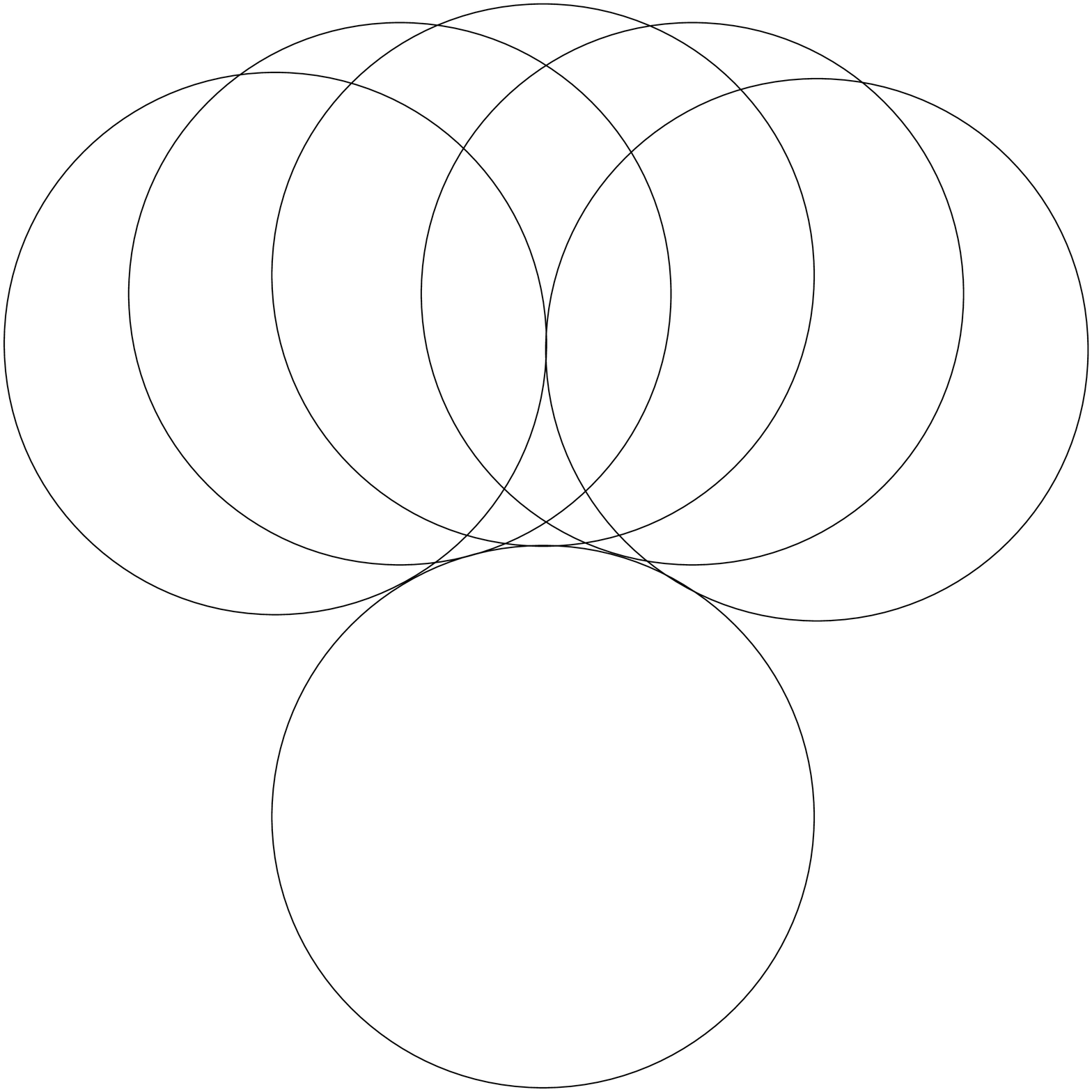}}
	\caption{\label{fig:sharp-uc} Complete arrangement of six unit circles with $v_0=4n-6$. 
	Points that look like touching points
	should be two intersection points between the respective circles. 
	Further circles can easily be added to obtain arrangements that meet the bound for arbitrary $n$. }
\end{figure}

We conjecture that Theorems \ref{thm:no4alpha}, \ref{thm:f0c}, and \ref{thm:k}
also hold if we drop the condition that the arrangement is $\beta$-free, i.e.,
for the improved bounds to hold it is sufficient that a complete arrangement is 
$\alpha^4$-free. However, the topology of these arrangements quickly becomes
rather involved so that we haven't yet succeeded in proving this.\footnote{The ``proof''
given in a prior version of this paper unfortunately turned out to contain a serious error.}
 As an $\alpha^4$-arrangement cannot be realized with unit circles, a proof of our 
conjecture would also imply that Theorems \ref{thm:no4alpha} and \ref{thm:k} 
especially hold for complete arrangements of unit circles. The improved bound
on $v_0$ would be sharp in this case (cf. Figure \ref{fig:sharp-uc}).

\bigskip

\end{document}